\documentclass[12pt]{article}
\usepackage{mathrsfs}
\usepackage{amssymb} \textwidth 150mm \textheight 220mm \oddsidemargin
15pt \evensidemargin 0pt \topmargin 0cm \headsep 0.3cm

\usepackage{amsmath}
\usepackage{graphicx}
\usepackage{indentfirst}
\usepackage{cite}
\usepackage{float}
\usepackage{CJK}

\begin{document}
\title{\Large\bf{Limit cycles appearing from the perturbation of differential systems with multiple switching curves
\thanks{E-mail address: jihua1113@163.com.} }}
\author{{Jihua Yang}\\
{\small \it School of Mathematics and Computer Science, Ningxia Normal University,}\\
 {\small \it  Guyuan, 756000, PR China}\\
 %{\small \it School of Mathematical Sciences, Beijing Normal University, }\\
%{\small \it Laboratory of Mathematics and Complex Systems, Ministry of Education,}\\
%{\small \it Beijing, 100875, PR China}
 }
\date{}
\maketitle \baselineskip=0.9\normalbaselineskip \vspace{-3pt}
\noindent
{\bf Abstract}\, This paper deals with the problem of limit cycle bifurcations for a piecewise near-Hamilton system with four regions separated by algebraic curves $y=\pm x^2$. By analyzing the obtained first order Melnikov function, we give an upper bound of the number of limit cycles which bifurcate from the period annulus around the origin under $n$-th degree polynomial perturbations. In the case $n=1$, we obtain that at least 4 (resp. 3) limit cycles can bifurcate from the period annulus if the switching curves are $y=\pm x^2$ (resp. $y=x^2$ or $y=-x^2$). The results also show that the number of switching curves affects the number of limit cycles.
\vskip 0.2 true cm
\noindent
{\bf Keywords}\, limit cycle; switching curve; Melnikov function

 \section{Introduction and main results}
 \setcounter{equation}{0}
\renewcommand\theequation{1.\arabic{equation}}

There are amount of non-smooth dynamical systems in nature and engineering areas where their complex  dynamic cannot be researched deeply in terms of classical smooth system theories due to the non-smooth dynamical systems' strong nonlinearities on the discontinuity sets (also called switching manifolds) \cite{BBC,K,C}. Therefore, it is very important to study dynamical behaviors, especially the bifurcations of limit cycle, of the non-smooth systems.

Piecewise smooth differential system is a kind of important non-smooth system which is based on non-smooth model.
% and can be expressed as
%\begin{eqnarray}\Big(\frac{dx}{dt},\frac{dy}{dt}\Big)=\begin{cases}
%\big(P^+(x,y),\  Q^+(x,y)\big),\quad C(x,y)>0,\\
%\big(P^-(x,y),\  Q^-(x,y)\big),\quad C(x,y)<0,\\
%\end{cases}\end{eqnarray}
%where $C(x,y)=0$ is called switching curve.
In the past few decades, many authors have been devoted to study the number of limit cycles of piecewise smooth differential systems with two zones separated by a switching line, see \cite{LH,LHR,X,LL,HZ,YZ17,LZ,CLT,GYC,CLYZ,SYZ} and the references quoted therein. The ways used in the aforementioned works are Melnikov function established in \cite{LH,HS} and averaging method developed in \cite{LMN,LNT}. Recently, Yang and Zhao \cite{YZ} used the Picard-Fuchs equation to calculate the first order Melnikov function of a kind of piecewise smooth differential systems with a switching line, which can reduce a lot of calculation work.

When the piecewise smooth differential systems are separated by finite many straight lines, there are some valuable results. Hu and Du \cite{HD} derived the first order Melnikov function for perturbed piecewise smooth differential systems with $m$ switching straight lines which can be used to study the number of limit cycles for these systems. Xiong \cite{X18} investigated the limit cycle bifurcation in perturbations of piecewise smooth Hamiltonian systems with switching lines $x=0$ and $y=0$ via multiple parameters. By using the averaging method of first order, Itikawa et al. \cite{IL} obtained  the upper bounds of the number of limit cycles bifurcating from the periodic orbits of two kind of isochronous systems, when they are perturbed inside the discontinuous quadratic and cubic polynomials differential systems, respectively. Akhmet and Aru\u{g}aslan \cite{AA} generalized the problem of Hopf bifurcation for a planar non-smooth system by considering discontinuities on finitely many nonlinear curves emanating from a vertex. For more results, one can see \cite{SD,LYD,LPP,WHC,DL,LT} and the references therein.

In the present paper, motivated by the above references, we will study the number of limit cycles for Hamilton system under perturbations of piecewise polynomials of degree $n$ with two switching curves $y=\pm x^2$. More precisely, we consider the following perturbed piecewise smooth differential system with four zones
\begin{eqnarray}
\left(  \begin{array}{c}
          \dot{x} \\          \dot{y}
          \end{array} \right)
=\begin{cases}
 \left(
  \begin{array}{c}
          y+\varepsilon f_1(x,y) \\
          -x+\varepsilon g_1(x,y)
          \end{array} \right), \ x^2>y>-x^2,\ x>0,\\[0.6truecm]
  \left(  \begin{array}{c}
         y+\varepsilon f_2(x,y) \\
          -x+\varepsilon g_2(x,y)
           \end{array}
 \right),\ \  y<-x^2,\\[0.6truecm]
  \left(  \begin{array}{c}
         y+\varepsilon f_3(x,y) \\
          -x+\varepsilon g_3(x,y)
           \end{array}
 \right),\ \ x^2>y>-x^2,\ x<0,\\[0.6truecm]
 \left(  \begin{array}{c}
         y+\varepsilon f_4(x,y) \\
          -x+\varepsilon g_4(x,y)
           \end{array}
 \right),\ y>x^2,
  \end{cases}
 \end{eqnarray}
where $$f_k(x,y)=\sum\limits_{i+j=0}^na^k_{i,j}x^iy^j,\ \ g_k(x,y)=\sum\limits_{i+j=0}^nb^k_{i,j}x^iy^j,\ k=1,2,3,4.$$
The first integral of system (1.1) for ${\varepsilon=0}$ is
\begin{eqnarray}
\begin{aligned}
&H^1(x,y)=\frac{1}{2}(x^2+y^2)=\frac{h}{2},\ \ x^2>y>-x^2,\ x>0,\\
&H^2(x,y)=\frac{1}{2}(x^2+y^2)=\frac{h}{2},\ \ y<-x^2,\\
&H^3(x,y)=\frac{1}{2}(x^2+y^2)=\frac{h}{2},\ \ x^2>y>-x^2,\ x<0,\\
&H^4(x,y)=\frac{1}{2}(x^2+y^2)=\frac{h}{2},\ \ y>x^2\\
\end{aligned}
\end{eqnarray}
with $h\in(0,+\infty)$. When $\varepsilon=0$, system (1.1) has a family of  periodic orbits as follows
$$\begin{aligned}
L_h=&\big\{(x,y)|H^1(x,y)=\frac{h}{2},\ x^2>y>-x^2,\ x>0\big\}\\
&\cup\big\{(x,y)|H^2(x,y)=\frac{h}{2},\ y<-x^2\big\}\\
&\cup\big\{(x,y)|H^3(x,y)=\frac{h}{2},\ x^2>y>-x^2,\ x<0\big\}\\
&\cup\big\{(x,y)|H^4(x,y)=\frac{h}{2},\ y>x^2\big\}\\
:=&L^1_h\cup L^2_h\cup L^3_h\cup L^4_h
\end{aligned}$$
with $h\in(0,+\infty)$, see Fig\,.1.
%\begin{figure}[htbp]
% \centering \includegraphics[width=3in]{1.eps}
%\begin{center}
%{\small{\bf Fig.\,1.} The phase portrait of system (1.1) for ${\varepsilon=0}$. }
%\end{center}
%\end{figure}

 Using methods of Theorem 1.1 in \cite{LH} and Theorem 4.1 in \cite{HD}, we can easily obtain the first order Melnikov function of system (1.1) which can be described in the following Proposition.
 \vskip 0.2 true cm

\noindent
{\bf Proposition 1.1.} {\it The first order Melnikov function of system (1.1) is
\begin{eqnarray}
\begin{aligned}
M(h)=&\Phi_1(h)\int_{\widehat{AB}}g_1(x,y)dx-f_1(x,y)dy+\Phi_2(h)\int_{\widehat{BC}}g_2(x,y)dx-f_2(x,y)dy\\
&+\Phi_3(h)\int_{\widehat{CD}}g_3(x,y)dx-f_3(x,y)dy+\Phi_4(h)\int_{\widehat{DA}}g_4(x,y)dx-f_4(x,y)dy,
\end{aligned}
\end{eqnarray}
where
$$\begin{aligned}
\Phi_1(h)=&\frac{\big[H^1_x(A)+H^1_y(A)f'_+(a(h))\big]\big[H^2_x(B)+H^2_y(B)f'_-(b(h))\big]}
{\big[H^4_x(A)+H^4_y(A)f'_+(a(h))\big]\big[H^1_x(B)+H^1_y(B)f'_-(b(h))\big]}\\
&\times \frac{\big[H^3_x(C)+H^3_y(C)f'_-(c(h))\big]\big[H^4_x(D)+H^4_y(D)f'_+(d(h))\big]}{\big[H^2_x(C)+H^2_y(C)f'_-(c(h))\big]
\big[H^3_x(D)+H^3_y(D)f'_+(d(h))\big]},\\
\Phi_2(h)=&\frac{\big[H^1_x(A)+H^1_y(A)f'_+(a(h))\big]\big[H^3_x(C)+H^3_y(C)f'_-(c(h))\big]}
{\big[H^4_x(A)+H^4_y(A)f'_+(a(h))\big]\big[H^2_x(C)+H^2_y(C)f'_-(c(h))\big]}\\
&\times \frac{H^4_x(D)+H^4_y(D)f'_+(d(h))}{H^3_x(D)+H^3_y(D)f'_+(d(h))},\\
\Phi_3(h)=&\frac{\big[H^1_x(A)+H^1_y(A)f'_+(a(h))\big]\big[H^4_x(D)+H^4_y(D)f'_+(d(h))\big]}
{\big[H^4_x(A)+H^4_y(A)f'_+(a(h))\big]\big[H^3_x(D)+H^3_y(D)f'_+(d(h))\big]},\\
\Phi_4(h)=&\frac{H^1_x(A)+H^1_y(A)f'_+(a(h))}{H^4_x(A)+H^4_y(A)f'_+(a(h))},\\
\end{aligned}$$
$f_\pm(x)=\pm x^2$ and $a(h)$, $b(h)$, $c(h)$ and $d(h)$ are the abscissas of points $A$, $B$, $C$ and $D$, respectively.
Further, if $M(h_0)=0$ and $M'(h_0)\neq0$ for some $h_0\in(0,+\infty)$, then for $|\varepsilon|$ small enough (1.1) has a unique limit cycle near $L_{h_0}$. If $h_0$ is a zero of $M(h)$ having an odd multiplicity, then for $|\varepsilon|$ small enough (1.1) has at least one limit cycle near $L_{h_0}$. Also, if $M(h)$ has at most $k$ zeros counting multiplicity in $h$ on the interval $(0,+\infty)$, then system (1.1) has at most $k$ limit cycles bifurcating from the annulus $\bigcup\limits_{h\in(0,+\infty)}L_h$.}
 \vskip 0.2 true cm

Let $H(n)$ denotes the upper bound of the number of limit cycles bifurcating from the period annulus around the origin for all possible polynomials $f_l(x,y)$ and $g_l(x,y)$ up to the first order Melnikov function, taking into account the multiplicity. Our main results are the following two theorems.
 \vskip 0.2 true cm

\noindent
{\bf Theorem 1.1.} {\it If the switching curves are $y=\pm x^2$, then
$$H(1)=4;\ \ H(n)\leq 2n+5[\frac{n-1}{2}]+4\ for\ n\geq2.$$}

\noindent
{\bf Theorem 1.2.} {\it If the switching curve is $y=x^2$ or $y=-x^2$, then
$$H(1)=3;\  \ H(n)\leq 2n+5[\frac{n-1}{2}]+4\ for\ n\geq2.$$
Notice that the equal sign ``=" means that the upper bound is reached.}
\vskip 0.2 true cm
\noindent
{\bf Remark 1.1.} If the switching curves are $m\geq2$ straight lines intersecting at the origin (see Fig.\,2), Carvalho, Llibre and Tonon \cite{CLT} proved that the maximum number of limit cycles bifurcating from the period annulus around the origin is $n$. Comparing with the results in Theorems 1.1 and 1.2, we know that the shapes of the switching curves have an essential effect on the number of limit cycles.
%\begin{figure}[htbp]
% \centering \includegraphics[width=3in]{3.eps}
%\begin{center}
%{\small{\bf Fig.\,2.} The phase portrait of system (1.1)$|_{\varepsilon=0}$ with $m$ switching lines. }
%\end{center}
%\end{figure}
\vskip 0.2 true cm

The layout of the rest of the paper is as follows. We first obtain the algebraic structure of the first order Melnikov function $M(h)$ in section 2. The main results will be proved in sections 3 and 4.

\section{Algebraic structure of $M(h)$}
 \setcounter{equation}{0}
\renewcommand\theequation{2.\arabic{equation}}

In the following, we will obtain the algebraic structure of $M(h)$ of system (1.1). By straightforward calculation, we have
$\Phi_i(h)=1,\ i=1,2,3,4.$
%Noting that the coordinates of $A_i$ for $i=1,2,3,4,5,6$ are
%$$\begin{aligned}&A_{1,6}=(\frac{3-\sqrt{5+4h}}{2},\pm\sqrt{\frac{3-\sqrt{5+4h}}{2}}),\ A_{2,5}=(0,\pm\sqrt{1-h}),\\ &A_{3,4}=(\frac{1-\sqrt{5-4h}}{2},\pm\sqrt{\frac{\sqrt{5-4h}-1}{2}}),\\
%&A_4=(\frac{1-\sqrt{1+4(1-h)}}{2},-\sqrt{\frac{\sqrt{1+4(1-h)}-1}{2}}),\\
%&A_6=(\frac{3-\sqrt{9-4(1-h)}}{2},-\sqrt{\frac{3-\sqrt{9-4(1-h)}}{2}}),\\
%\end{aligned}$$
%Hence, we have
Hence, the first Melnikov function $M(h)$ of system (1.1) has the form
\begin{eqnarray}
\begin{aligned}
M(h)=&\int_{\widehat{AB}}g_1(x,y)dx-f_1(x,y)dy+\int_{\widehat{BC}}g_2(x,y)dx-f_2(x,y)dy\\
&+\int_{\widehat{CD}}g_3(x,y)dx-f_3(x,y)dy+\int_{\widehat{DA}}g_4(x,y)dx-f_4(x,y)dy.
\end{aligned}
\end{eqnarray}
  For $h\in(0,+\infty)$ and $i,j\in\mathbb{N}$, we denote
$$\begin{aligned}
&I_{i,j}(h)=\int_{\widehat{AB}}x^iy^jdx,\ \ J_{i,j}(h)=\int_{\widehat{BC}}x^iy^jdx,\\
&\tilde{I}_{i,j}(h)=\int_{\widehat{CD}}x^iy^jdx,\ \ \tilde{J}_{i,j}(h)=\int_{\widehat{DA}}x^iy^jdx.
\end{aligned}$$
\vskip 0.2 true cm

\noindent
{\bf Lemma 2.1.}\, {\it The first order Melnikov function $M(h)$ can be written as
\begin{eqnarray}\begin{aligned}
M(h)=&\sum\limits_{i+j=0}^n\tau_{i,j}I_{i,j}(h)+\sum\limits_{i+j=0}^n\sigma_{i,j}J_{i,j}(h)+
\tilde{\phi}_{2n+\frac{3+(-1)^n}{2}}\Big(\sqrt{\frac{\sqrt{1+4h}-1}{2}}\Big),
\end{aligned}\end{eqnarray}
where $\tau_{i,j}$ and $\sigma_{i,j}$ are arbitrary constants, $\phi_{l}(u)$ is a polynomial of $u$ of degree at most $l$.}
\vskip 0.2 true cm

\noindent
{\bf Proof.}  Let $\Omega$ be the interior of $\widehat{AB}\cup \widehat{BO}\cup \widehat{OA}$, see Fig.\,1. Using the Green's Formula, we have
\begin{eqnarray*}
\begin{aligned}
\int_{\widehat{AB}}x^iy^jdx
&=\oint_{\widehat{AB}\cup \widehat{BO}\cup \widehat{OA}}x^iy^jdx-\int_{\widehat{BO}}x^iy^jdx-\int_{\widehat{OA}}x^iy^jdx\\
&=j\iint\limits_{\Omega} x^iy^{j-1}dxdy-\frac{1-(-1)^j}{i+2j+1}\Big(\frac{\sqrt{1+4h}-1}{2}\Big)^\frac{i+2j+1}{2},\\
\int_{\widehat{AB}}x^iy^jdy
&=\oint_{\widehat{AB}\cup \widehat{BO}\cup \widehat{OA}}x^iy^jdy-\int_{\widehat{BO}}x^iy^jdy-\int_{\widehat{OA}}x^iy^jdy\\
&=-i\iint\limits_{\Omega} x^{i-1}y^{j}dxdy-\frac{2[1+(-1)^j]}{i+2j+2}\Big(\frac{\sqrt{1+4h}-1}{2}\Big)^\frac{i+2j+2}{2},\\
\end{aligned}
\end{eqnarray*}
which imply that
\begin{eqnarray}
\int_{\widehat{AB}}x^iy^jdy=-\frac{i}{j+1}I_{i-1,j+1}(h)-\frac{1+(-1)^j}{j+1}\Big(\frac{\sqrt{1+4h}-1}{2}\Big)^\frac{i+2j+2}{2}.
\end{eqnarray}
In a similar way, we have
\begin{eqnarray}
\begin{aligned}
&\int_{\widehat{BC}}x^iy^jdy=-\frac{i}{j+1}J_{i-1,j+1}(h)+(-1)^j\frac{1-(-1)^i}{j+1}\Big(\frac{\sqrt{1+4h}-1}{2}\Big)^\frac{i+2j+2}{2},\\
&\int_{\widehat{CD}}x^iy^jdy=-\frac{i}{j+1}\tilde{I}_{i-1,j+1}(h)+(-1)^i\frac{1+(-1)^j}{j+1}\Big(\frac{\sqrt{1+4h}-1}{2}\Big)^\frac{i+2j+2}{2},\\
&\int_{\widehat{DA}}x^iy^jdy=-\frac{i}{j+1}\tilde{J}_{i-1,j+1}(h)+\frac{1-(-1)^i}{j+1}\Big(\frac{\sqrt{1+4h}-1}{2}\Big)^\frac{i+2j+2}{2}.
\end{aligned}
\end{eqnarray}

From (2.1), (2.3) and (2.4), we obtain
\begin{eqnarray*}
\begin{aligned}
M(h)=&\sum\limits_{i+j=0}^nb^1_{i,j}I_{i,j}(h)+\sum\limits_{i+j=0}^n\frac{i}{j+1}a^1_{i,j}I_{i-1,j+1}(h)\\&+
\sum\limits_{i+j=0}^n\frac{1+(-1)^j}{j+1}a^1_{i,j}\Big(\frac{\sqrt{1+4h}-1}{2}\Big)^\frac{i+2j+2}{2}\\
&+\sum\limits_{i+j=0}^nb^2_{i,j}J_{i,j}(h)+\sum\limits_{i+j=0}^n\frac{i}{j+1}a^2_{i,j}J_{i-1,j+1}(h)\\&+
\sum\limits_{i+j=0}^n(-1)^j\frac{(-1)^i-1}{j+1}a^2_{i,j}\Big(\frac{\sqrt{1+4h}-1}{2}\Big)^\frac{i+2j+2}{2}\\
&+\sum\limits_{i+j=0}^nb^3_{i,j}\tilde{I}_{i,j}(h)+\sum\limits_{i+j=0}^n\frac{i}{j+1}a^3_{i,j}\tilde{I}_{i-1,j+1}(h)\\&-
\sum\limits_{i+j=0}^n(-1)^i\frac{1+(-1)^j}{j+1}a^3_{i,j}\Big(\frac{\sqrt{1+4h}-1}{2}\Big)^\frac{i+2j+2}{2}\\
&+\sum\limits_{i+j=0}^nb^4_{i,j}\tilde{J}_{i,j}(h)+\sum\limits_{i+j=0}^n\frac{i}{j+1}a^4_{i,j}\tilde{J}_{i-1,j+1}(h)\\&+
\sum\limits_{i+j=0}^n\frac{(-1)^i-1}{j+1}a^4_{i,j}\Big(\frac{\sqrt{1+4h}-1}{2}\Big)^\frac{i+2j+2}{2}.
\end{aligned}
\end{eqnarray*}
It is easy to check that
\begin{eqnarray*}
\tilde{I}_{i,j}(h)=(-1)^{i+j+1}I_{i,j}(h),\ \ \tilde{J}_{i,j}(h)=(-1)^jJ_{i,j}(h).
\end{eqnarray*}
Hence,
\begin{eqnarray}
\begin{aligned}
M(h)=&\sum\limits_{i+j=0}^n\tau_{i,j}I_{i,j}(h)+\sum\limits_{i+j=0}^n\sigma_{i,j}J_{i,j}(h)
+\sum\limits_{i+j=0}^n\rho_{i,j}\Big(\frac{\sqrt{1+4h}-1}{2}\Big)^\frac{i+2j+2}{2}\\
:=&\sum\limits_{i+j=0}^n\tau_{i,j}I_{i,j}(h)+\sum\limits_{i+j=0}^n\sigma_{i,j}J_{i,j}(h)
+\tilde{\phi}_{2n+\frac{3+(-1)^n}{2}}\Big(\sqrt{\frac{\sqrt{1+4h}-1}{2}}\Big),\\
\end{aligned}
\end{eqnarray}
where $\tau_{i,j}$, $\sigma_{i,j}$ and $\rho_{i,j}$ are arbitrary constants and $\phi_{l}(u)$ is a polynomial of $u$ of degree at most $l$.
\vskip 0.2 true cm

\noindent
{\bf Lemma 2.2.}\, {\it For $h\in(0,+\infty)$ and $l+m\geq1$, we have}
\vskip 0.2 true cm

\noindent
(i) {\it If $n=2l+2m$, then
\begin{eqnarray}
\begin{aligned}
&J_{2l,2m}(h)=\tilde{\gamma}_{l,m}h^{[\frac{n}{2}]}J_{0,0}(h)+\sum\limits_{k=2}^{n}\hat{\varphi}_{[\frac{n-k}{2}]}(h)\Big(\frac{\sqrt{1+4h}-1}{2}\Big)^{k+\frac{1}{2}}.\\
\end{aligned}
\end{eqnarray}}

\noindent
(ii) {\it If $n=2l+2m+1$, then
\begin{eqnarray}
\begin{aligned}
&I_{2l,2m+1}(h)=\tilde{\alpha}_{l,m}h^{[\frac{n-1}{2}]}I_{0,1}(h)+\sum\limits_{k=3}^{n}\tilde{\varphi}_{[\frac{n-k}{2}]}(h)\Big(\frac{\sqrt{1+4h}-1}{2}\Big)^{k+\frac{1}{2}},\\
&J_{2l,2m+1}(h)=\tilde{\delta}_{l,m}h^{[\frac{n-1}{2}]}J_{0,1}(h)+\sum\limits_{k=3}^{n}\hat{\psi}_{[\frac{n-k}{2}]}(h)\Big(\frac{\sqrt{1+4h}-1}{2}\Big)^{k+\frac{1}{2}}.\\
\end{aligned}
\end{eqnarray}}

\noindent
(iii) {\it If $n=2l+2m+2$, then
\begin{eqnarray}
\begin{aligned}
&I_{2l+1,2m+1}(h)=\tilde{\beta}_{l,m}h^{[\frac{n-2}{2}]}I_{1,1}(h)+\sum\limits_{k=4}^{n}\tilde{\psi}_{[\frac{n-k}{2}]}(h)\Big(\frac{\sqrt{1+4h}-1}{2}\Big)^{k}.\\
\end{aligned}
\end{eqnarray}
where $[p]$ denotes the integer part of $p$, $\tilde{\alpha}_{l,m}$, $\tilde{\beta}_{l,m}$, $\tilde{\gamma}_{l,m}$ and $\tilde{\delta}_{l,m}$ are arbitrary constants, and $\tilde{\varphi}_{l}(h)$, $\hat{\varphi}_{l}(h)$, $\tilde{\psi}_{l}(h)$ and $\hat{\psi}_{l}(h)$ are polynomials of $h$ with degree no more than $l$.}
\vskip 0.2 true cm
\noindent
{\bf Proof.}  Since the integral path $\widehat{AB}$ (resp. $\widehat{BC}$) are symmetrical with respect to $x$-axis (resp. $y$-axis), $I_{i,2l}(h)=0$, $J_{2l+1,j}(h)=0$.

Without loss of generality, we only prove the first equality in (2.7), and the others can be shown in a similar way.
Differentiating $H^1(x,y)=\frac{h}{2}$ defined in (1.2) with respect to $x$, we have
\begin{eqnarray}
x+y\frac{\partial y}{\partial x}=0.
\end{eqnarray}
Multiplying (2.9) by $x^{i-1}y^{j}dx$, integrating over $\widehat{AB}$ and noting that (2.3), we have
\begin{eqnarray}\begin{aligned}
I_{i,j}(h)=&\frac{i-1}{j+2}I_{i-2,j+2}(h)+\frac{1-(-1)^j}{j+2}\Big(\frac{\sqrt{1+4h}-1}{2}\Big)^\frac{i+2j+3}{2}.
\end{aligned}\end{eqnarray}

On the other hand, multiplying $H^1(x,y)=\frac{h}{2}$ defined in (1.2) by $x^{i}y^{j-2}dy$ and integrating over $\widehat{AB}$ yield
\begin{eqnarray}
I_{i,j}(h)=hI_{i,j-2}(h)-I_{i+2,j-2}(h).
\end{eqnarray}
From (2.10) and (2.11), we have
\begin{eqnarray}\begin{aligned}
I_{i,j}(h)=\frac{1}{i+j+1}\Big[(i-1)hI_{i-2,j}(h)+\big(1-(-1)^j\big)\Big(\frac{\sqrt{1+4h}-1}{2}\Big)^\frac{i+2j+3}{2}\Big]
\end{aligned}\end{eqnarray}
and
\begin{eqnarray}\begin{aligned}
I_{i,j}(h)=\frac{1}{i+j+1}\Big[jhI_{i,j-2}(h)-\big(1-(-1)^j\big)\Big(\frac{\sqrt{1+4h}-1}{2}\Big)^\frac{i+2j+1}{2}\Big].
\end{aligned}\end{eqnarray}

We will prove the conclusion by induction on $l+m=p$. It could be noticed that $p=1$ corresponds to $(i,j)=(0,3)$ and $(2,1)$ and $p=2$ corresponds to $(i,j)=(0,5)$, $(2,3)$ and $(4,1)$. Hence, in view of (2.12) and (2.13), we have
\begin{eqnarray}
\begin{cases}
I_{0,3}(h)=\frac{3}{4}hI_{0,1}(h)-\frac{1}{2}\Big(\frac{\sqrt{1+4h}-1}{2}\Big)^\frac{7}{2},\\
I_{2,1}(h)=\frac{1}{4}hI_{0,1}(h)+\frac{1}{2}\Big(\frac{\sqrt{1+4h}-1}{2}\Big)^\frac{7}{2},\\
I_{0,5}(h)=\frac{5}{8}h^2I_{0,1}(h)-\frac{5}{12}h\Big(\frac{\sqrt{1+4h}-1}{2}\Big)^\frac{7}{2}
-\frac{1}{3}\Big(\frac{\sqrt{1+4h}-1}{2}\Big)^\frac{11}{2},\\
I_{2,3}(h)=\frac{1}{8}h^2I_{0,1}(h)-\frac{1}{4}h\Big(\frac{\sqrt{1+4h}-1}{2}\Big)^\frac{7}{2}
-\frac{1}{3}\Big(\frac{\sqrt{1+4h}-1}{2}\Big)^\frac{9}{2},\\
I_{4,1}(h)=\frac{1}{8}h^2I_{0,1}(h)-\frac{1}{4}h\Big(\frac{\sqrt{1+4h}-1}{2}\Big)^\frac{7}{2}
+\frac{1}{3}\Big(\frac{\sqrt{1+4h}-1}{2}\Big)^\frac{9}{2}.
\end{cases}
\end{eqnarray}
which yield the conclusion for $p=1,2$. Now assume that the result holds for all $l+m\leq p-1$ ($p\geq2$). Then, for $l+m=p$, taking $(i,j)=(0,2p+1),(2,2p-1),\cdots,(2p-4,5),(2p-2,3)$ in (2.13) and $(i,j)=(2p,1)$ in (2.12), respectively, we obtain
\begin{eqnarray}
\left(\begin{matrix}
                 I_{0,2p+1}(h)\\
                 I_{2,2p-1}(h)\\
                                  \vdots\\
                  I_{2p-4,5}(h)\\
                  I_{2p-2,3}(h)\\
                  I_{2p,1}(h)\\
\end{matrix}\right)\ \
=\frac{1}{2p+2}\left(\begin{matrix}
                (2p+1)hI_{0,2p-1}(h)-2\Big(\frac{\sqrt{1+4h}-1}{2}\Big)^{2p+\frac{3}{2}}\\
                (2p-1)hI_{2,2p-3}(h)-2\Big(\frac{\sqrt{1+4h}-1}{2}\Big)^{2p+\frac{1}{2}}\\
                % (2k-3)hI_{4,2k-5}(h)-2\Big(\frac{\sqrt{1+4h}-1}{2}\Big)^{2k-\frac{1}{2}}\\
                 \vdots\\
5hI_{2p-4,3}(h)-2\Big(\frac{\sqrt{1+4h}-1}{2}\Big)^{p+\frac{7}{2}}\\
3hI_{2p-2,1}(h)-2\Big(\frac{\sqrt{1+4h}-1}{2}\Big)^{p+\frac{5}{2}}\\
(2p-1)hI_{2p-2,1}(h)+2\Big(\frac{\sqrt{1+4h}-1}{2}\Big)^{p+\frac{5}{2}}
\end{matrix}\right).
\end{eqnarray}

By inductive hypothesis and (2.15), we have for $l+m=p$
\begin{eqnarray*}
\begin{aligned}
I_{2l,2m+1}(h)=&h\Big[\tilde{\alpha}_{l,m}h^{[\frac{2p-1-1}{2}]}I_{0,1}(h)+\sum\limits_{k=3}^{2p-1}\tilde{\varphi}_{[\frac{2p-1-k}{2}]}(h)
\Big(\frac{\sqrt{1+4h}-1}{2}\Big)^{k+\frac{1}{2}}\Big]\\&-2\Big(\frac{\sqrt{1+4h}-1}{2}\Big)^{2p+\frac{3}{2}}\\
:=&\tilde{\alpha}_{l,m}h^{[\frac{2p}{2}]}I_{0,1}(h)+\sum\limits_{k=3}^{2p+1}\tilde{\varphi}_{[\frac{2p+1-k}{2}]}(h)
\Big(\frac{\sqrt{1+4h}-1}{2}\Big)^{k+\frac{1}{2}},
\end{aligned}
\end{eqnarray*}
where $\tilde{\psi}_{l}(h)$ is a polynomial of $h$ with degree no more than $l$. The proof is completed.\quad $\lozenge$
\vskip 0.2 true cm

Substituting (2.6)-(2.8) into (2.5), we get the algebraic structure of the first order Melnikov function $M(h)$.
\vskip 0.2 true cm

\noindent
{\bf Lemma 2.3.}\, {\it For $h\in(0,+\infty)$, the first order Melnikov function $M(h)$ can be written as
\begin{eqnarray}
\begin{aligned}
M(h)=&\alpha(h)I_{0,1}(h)+\beta(h)I_{1,1}(h)+\gamma(h)J_{0,0}(h)+\delta(h)J_{0,1}(h)\\&
+\phi_{2n+\frac{3+(-1)^n}{2}}\Big(\sqrt{\frac{\sqrt{1+4h}-1}{2}}\Big),
\end{aligned}
\end{eqnarray}
where $\alpha(h)$, $\beta(h)$, $\gamma(h)$ and $\delta(h)$ are polynomials of $h$ satisfying
$$\deg \alpha(h),\deg \delta(h)\leq [\frac{n-1}{2}],\ \deg \beta(h)\leq[\frac{n-2}{2}],\ \deg \gamma(h)\leq[\frac{n}{2}].$$}

 \section{Proof of the Theorem 1.1}
 \setcounter{equation}{0}
\renewcommand\theequation{3.\arabic{equation}}

In the following, we denote by $P_k(u)$,  $Q_k(u)$, $R_k(u)$ and $S_k(u)$ polynomials of $u$ with degree at most $k$.

By some straightforward calculations, we have
\begin{eqnarray}
\begin{aligned}
&J_{0,0}(h)=-2\sqrt{\frac{\sqrt{1+4h}-1}{2}},\\
&I_{1,1}(h)=\frac{2}{3}\Big(h-{\frac{\sqrt{1+4h}-1}{2}}\Big)^\frac{3}{2},\\
&I_{0,1}(h)=2\int_0^{\sqrt{h}}\sqrt{h-x^2}dx-2\int_0^{\sqrt{\frac{\sqrt{1+4h}-1}{2}}}\sqrt{h-x^2}dx,\\
&J_{0,1}(h)=2\int_0^{\sqrt{\frac{\sqrt{1+4h}-1}{2}}}\sqrt{h-x^2}dx.
\end{aligned}
\end{eqnarray}
 If $n\geq2$, substituting (3.1) into (2.16), we have
\begin{eqnarray}
\begin{aligned}
M(h)=&P_{[\frac{n-1}{2}]}(h)\int_0^{\sqrt{h}}\sqrt{h-x^2}dx+Q_{[\frac{n-1}{2}]}(h)\int_0^{\sqrt{\frac{\sqrt{1+4h}-1}{2}}}\sqrt{h-x^2}dx\\&
+R_{[\frac{n-2}{2}]}(h)\Big(h-{\frac{\sqrt{1+4h}-1}{2}}\Big)^\frac{3}{2}+S_{[\frac{n}{2}]}(h)\sqrt{\frac{\sqrt{1+4h}-1}{2}}\\
&+\phi_{2n+\frac{3+(-1)^n}{2}}\Big(\sqrt{\frac{\sqrt{1+4h}-1}{2}}\Big).
\end{aligned}
\end{eqnarray}
Let $x=\sqrt{h}t$, then $M(h)$ in (3.2) becomes
\begin{eqnarray}
\begin{aligned}
M(h)=&hP_{[\frac{n-1}{2}]}(h)\int_0^1\sqrt{1-t^2}dt+hQ_{[\frac{n-1}{2}]}(h)\int_0^{\sqrt{\frac{\sqrt{1+4h}-1}{2h}}}\sqrt{1-t^2}dt\\&
+R_{[\frac{n-2}{2}]}(h)\Big(h-{\frac{\sqrt{1+4h}-1}{2}}\Big)^\frac{3}{2}+S_{[\frac{n}{2}]}(h)\sqrt{\frac{\sqrt{1+4h}-1}{2}}\\
&+\phi_{2n+\frac{3+(-1)^n}{2}}\Big(\sqrt{\frac{\sqrt{1+4h}-1}{2}}\Big),
\end{aligned}
\end{eqnarray}
where $\int_0^1\sqrt{1-t^2}dt=\frac{\pi}{4}$.

Let $\sqrt{\frac{\sqrt{1+4h}-1}{2}}=u$, that is, $h=u^4+u^2$.  Then $M(h)$ in (3.3) can be written as
\begin{eqnarray}
\begin{aligned}
M(u)=&uP_{2n+1}(u)+(u^4+u^2)Q_{[\frac{n-1}{2}]}(u^4+u^2)\int_0^{\frac{1}{\sqrt{1+u^2}}}\sqrt{1-t^2}dt.
\end{aligned}
\end{eqnarray}
%where $P_{2n+2}(u)$ is a polynomial of $u$ with degree at most $2n+2$ and does not have constant term.
It is easy to check that $M(h)$ and $M(u)$ have the same number of zeros in $(0,+\infty)$. Suppose that $\Sigma_1=(0,+\infty)\setminus \{u\in(0,+\infty)|(u^4+u^2)Q_{[\frac{n-1}{2}]}(u^4+u^2)=0\}$. By direct computation, we obtain for $u\in\Sigma_1$
\begin{eqnarray}
\begin{aligned}
\frac{d}{du}\Big(\frac{M(u)}{(u^4+u^2)Q_{[\frac{n-1}{2}]}(u^4+u^2)}\Big)
=&\Big(\frac{P_{2n+1}(u)}{(u^3+u)Q_{[\frac{n-1}{2}]}(u^4+u^2)}\Big)'-\frac{u^2}{(1+u^2)^2}\\[0.5truecm]
=&\frac{P_{2n+4[\frac{n-1}{2}]+3}(u)}{(u^3+u)^2Q^2_{[\frac{n-1}{2}]}(u^4+u^2)}.\\
\end{aligned}
\end{eqnarray}
Therefore, by Rolle's theorem,, $M(u)$ has at most $2n+5[\frac{n-1}{2}]+4$ zeros in $(0,+\infty)$, so does $M(h)$. That is,
$$H(n)\leq 2n+5[\frac{n-1}{2}]+4,\ n\geq2.$$

 If $n=1$, we have
\begin{eqnarray}\begin{aligned}
M(h)=&(b^1_{0,1}+b^3_{0,1}+a^1_{1,0}+a^3_{1,0})I_{0,1}(h)+(b^2_{0,0}+b^4_{0,0})J_{0,0}(h)\\
&+(b^2_{0,1}-b^4_{0,1}+a^2_{1,0}+a^4_{1,0})J_{0,1}(h)+2(a^1_{0,0}-a^3_{0,0})\frac{\sqrt{1+4h}-1}{2}\\
&+2(a^1_{1,0}-a^2_{1,0}+a^3_{1,0}-a^4_{1,0})\Big(\frac{\sqrt{1+4h}-1}{2}\Big)^\frac{3}{2},
\end{aligned}\end{eqnarray}
Let $\sqrt{\frac{\sqrt{1+4h}-1}{2}}=u$, we get
\begin{eqnarray}\begin{aligned}
M(u)=&\frac{\pi}{2}(b^1_{0,1}+b^3_{0,1}+a^1_{1,0}+a^3_{1,0})u^4+2(a^1_{1,0}-a^2_{1,0}+a^3_{1,0}-a^4_{1,0})u^3\\&+
[\frac{\pi}{2}(b^1_{0,1}+b^3_{0,1}+a^1_{1,0}+a^3_{1,0})+2(a^1_{0,0}-a^3_{0,0})]u^2-2(b^2_{0,0}+b^4_{0,0})u\\
&+2(b^2_{0,1}-b^4_{0,1}+a^2_{1,0}+a^4_{1,0}-b^1_{0,1}-b^3_{0,1}-a^1_{1,0}-a^3_{1,0})(u^4+u^2)\int_0^{\frac{1}{\sqrt{1+u^2}}}\sqrt{1-t^2}dt\\
:=&\lambda_4u^4+\lambda_3u^3+\lambda_2u^2+\lambda_1u+\lambda_0(u^4+u^2)\int_0^{\frac{1}{\sqrt{1+u^2}}}\sqrt{1-t^2}dt.
\end{aligned}\end{eqnarray}
We can prove that $M(u)$ in (3.7) has at most 4 zeros in $(0,+\infty)$ by using the same method above.

The Taylor expansion of function $\int_0^{\frac{1}{\sqrt{1+u^2}}}\sqrt{1-t^2}dt$ in the variable $u$, around $u=0$, is
\begin{eqnarray}\int_0^{\frac{1}{\sqrt{1+u^2}}}\sqrt{1-t^2}dt=\frac{\pi}{4}-\frac{1}{3}u^3+o(u^3).\end{eqnarray}
Therefore,
\begin{eqnarray}\begin{aligned}
M(u)=-\frac{\lambda_0}{3}u^5+(\lambda_4+\frac{\pi}{4}\lambda_0)u^4+\lambda_3u^3+(\lambda_2+\frac{\pi}{4}\lambda_0)u^2+\lambda_1 u+o(u^5).
\end{aligned}\end{eqnarray}
It is easy to check that the determinant of the Jacobian
 $$\det\frac{\partial(\lambda_1,\mu_2,\lambda_3,\mu_4,\mu_5)}{\partial(b_{0,1}^2,b_{0,0}^2,a_{0,0}^1,a_{1,0}^2,
b_{0,1}^1)}=\frac{8}{3}\pi,$$
where $\mu_2=\lambda_2+\frac{\pi}{4}\lambda_0,\mu_4=\lambda_4+\frac{\pi}{4}\lambda_0,\mu_5=-\frac{1}{3}\lambda_0.$
That is, $\lambda_1,\mu_2,\lambda_3,\mu_4$ and $\mu_5$ can be chosen arbitrarily. Hence, we can choose $\lambda_1,\mu_2,\lambda_3,\mu_4$ and $\mu_5$ appropriately such that $M(u)$ in (3.9) has 4 zeros in $(0,+\infty)$. This ends the proof of Theorem 1.1. \quad $\lozenge$

 \section{Proof of the Theorem 1.2}
 \setcounter{equation}{0}
\renewcommand\theequation{4.\arabic{equation}}

 Without loss of generality, we only consider the case that the switching curve is $y=x^2$. The other case can be shown similarly. Thus, system (1.1) can be written as
\begin{eqnarray}
\left(  \begin{array}{c}
          \dot{x} \\          \dot{y}
          \end{array} \right)
=\begin{cases}
 \left(
  \begin{array}{c}
          y+\varepsilon f_1(x,y) \\
          -x+\varepsilon g_1(x,y)
          \end{array} \right), \ y< x^2,\\[0.6truecm]
   \left(  \begin{array}{c}
         y+\varepsilon f_4(x,y) \\
          -x+\varepsilon g_4(x,y)
           \end{array}
 \right),\ y>x^2.
  \end{cases}
 \end{eqnarray}
Similar to Proposition 1.1, the first order Melnikov function of system (4.1) has the form
\begin{eqnarray}
M(h)=\int_{\widehat{AD}}g^1(x,y)dx-f^1(x,y)dy+\int_{\widehat{DA}}g^4(x,y)dx-f^4(x,y)dy,
\end{eqnarray}
and the number of zeros of the above Melnikov function for $|\varepsilon|$ small enough controls the number of limit cycles of system (4.1) bifurcating from the period annulus, see Fig.\,3.
%\begin{figure}[htbp]
% \centering \includegraphics[width=3in]{2.eps}
%\begin{center}
%{\small{\bf Fig.\,3.} The phase portrait of system (4.1)$|_{\varepsilon=0}$ with a switching curve $y=x^2$. }
%\end{center}
%\end{figure}

Similar to (2.3), we have
\begin{eqnarray}
\begin{aligned}
\int_{\widehat{AD}}x^iy^jdy=-\frac{i}{j+1}U_{i-1,j+1}(h)+\frac{(-1)^i-1}{j+1}\Big(\frac{\sqrt{1+4h}-1}{2}\Big)^\frac{i+2j+2}{2},\\
\int_{\widehat{DA}}x^iy^jdy=-\frac{i}{j+1}V_{i-1,j+1}(h)-\frac{(-1)^i-1}{j+1}\Big(\frac{\sqrt{1+4h}-1}{2}\Big)^\frac{i+2j+2}{2},
\end{aligned}
\end{eqnarray}
where $$U_{i,j}(h)=\int_{\widehat{AD}}x^iy^jdx,\ \ V_{i,j}(h)=\int_{\widehat{DA}}x^iy^jdx.$$
Since $\widehat{AD}$ and $\widehat{DA}$ are symmetric with respect to the $y$-axis,
$U_{2l+1,j}(h)=V_{2l+1,j}(h)=0$. Therefore, by (4.2) and (4.3), we have
\begin{eqnarray*}
\begin{aligned}
M(h)=&\sum\limits_{i+j=0}^nb^1_{i,j}U_{i,j}(h)+\sum\limits_{i+j=0}^n\frac{i}{j+1}a^1_{i,j}U_{i-1,j+1}(h)\\&-
\sum\limits_{i+j=0}^n\frac{(-1)^i-1}{j+1}a^1_{i,j}\Big(\frac{\sqrt{1+4h}-1}{2}\Big)^\frac{i+2j+2}{2}\\
&+\sum\limits_{i+j=0}^nb^4_{i,j}V_{i,j}(h)+\sum\limits_{i+j=0}^n\frac{i}{j+1}a^4_{i,j}V_{i-1,j+1}(h)\\&+
\sum\limits_{i+j=0}^n\frac{(-1)^i-1}{j+1}a^4_{i,j}\Big(\frac{\sqrt{1+4h}-1}{2}\Big)^\frac{i+2j+2}{2}\\
=&\sum\limits_{i+j=0}^n\bar{\tau}_{i,j}U_{i,j}(h)+\sum\limits_{i+j=0}^n\bar{\sigma}_{i,j}V_{i,j}(h)
+\Big(\frac{\sqrt{1+4h}-1}{2}\Big)^{\frac{3}{2}}\bar{\varphi}_{n-1}\Big(\frac{\sqrt{1+4h}-1}{2}\Big),
\end{aligned}
\end{eqnarray*}
where $\bar{\tau}_{i,j}$ and $\bar{\sigma}_{i,j}$ are arbitrary constants and $\bar{\varphi}_{n-1}(u)$ is polynomial of $u$ of degree no more that $n-1$.

Similar to the proof of Lemma 2.2, we can obtain the algebraic structure of $M(h)$.
\vskip 0.2 true cm

\noindent
{\bf Lemma 4.1.}\, {\it For $h\in(0,+\infty)$, the first order Melnikov function of system (4.1) can be written as
\begin{eqnarray}
\begin{aligned}
M(h)=&\alpha(h)U_{0,0}(h)+\beta(h)U_{0,1}(h)+\gamma(h)V_{0,0}(h)+\delta(h)V_{0,1}(h)\\&
+\Big(\frac{\sqrt{1+4h}-1}{2}\Big)^{\frac{3}{2}}{\varphi}_{n-1}\Big(\frac{\sqrt{1+4h}-1}{2}\Big),
\end{aligned}
\end{eqnarray}
where ${\varphi}_{n-1}(u)$ is polynomial of $u$ of degree no more that $n-1$ and $\alpha(h)$, $\beta(h)$, $\gamma(h)$ and $\delta(h)$ are polynomials of $h$ satisfying
$$\deg \alpha(h),\deg \gamma(h)\leq [\frac{n}{2}],\ \deg \beta(h),\deg\delta(h)\leq[\frac{n-1}{2}].$$}

\noindent
{\bf Proof of the Theorem 1.2.} By direct computation, we have
\begin{eqnarray}
\begin{aligned}
&U_{0,0}(h)=-2\sqrt{\frac{\sqrt{1+4h}-1}{2}},\  U_{0,1}(h)=4\int_0^{\sqrt{h}}\sqrt{h-x^2}dx-2\int_0^{\sqrt{\frac{\sqrt{1+4h}-1}{2}}}\sqrt{h-x^2}dx,\\
&V_{0,0}(h)=2\sqrt{\frac{\sqrt{1+4h}-1}{2}},\ \ \ \ V_{0,1}(h)=2\int_0^{\sqrt{\frac{\sqrt{1+4h}-1}{2}}}\sqrt{h-x^2}dx.
\end{aligned}
\end{eqnarray}
 Substituting (4.5) into (4.4) and let $\sqrt{\frac{\sqrt{1+4h}-1}{2}}=u$, we obtain
\begin{eqnarray*}
M(u)=uP_{2n+1}(u)+(u^4+u^2)Q_{[\frac{n-1}{2}]}(u^4+u^2)\int_0^{\frac{1}{\sqrt{1+u^2}}}\sqrt{1-t^2}dt,\ u\in(0,+\infty).
\end{eqnarray*}
Following the lines of the proof of Theorem 1.2, we obtain that $M(h)$ has at most $2n+5[\frac{n-1}{2}]+4$ zeros.

If $n=1$, then we have
$$\begin{aligned}M(h)=&b^1_{0,0}U_{0,0}(h)+(b^1_{0,1}+a^1_{1,0})U_{0,1}(h)+b_{0,0}^4V_{0,0}(h)+(b^4_{0,1}+a^4_{1,0})V_{0,1}(h)\\
&+2(a^1_{1,0}-a^4_{1,0})\Big(\frac{\sqrt{1+4h}-1}{2}\Big)^\frac{3}{2}\\
=&2(b_{0,0}^4-b_{0,0}^1)\sqrt{\frac{\sqrt{1+4h}-1}{2}}+2(a_{1,0}^1-a_{1,0}^4)\Big(\frac{\sqrt{1+4h}-1}{2}\Big)^\frac{3}{2}\\
&+4(b_{0,1}^1+a_{1,0}^1)\int_0^{\sqrt{h}}\sqrt{h-x^2}dx\\&
+2(b_{0,1}^4+a_{1,0}^4-b_{0,1}^1-a_{1,0}^1)\int_0^{\sqrt{\frac{\sqrt{1+4h}-1}{2}}}\sqrt{h-x^2}dx.\end{aligned}$$
Let $\sqrt{\frac{\sqrt{1+4h}-1}{2}}=u$, we obtain
\begin{eqnarray}\begin{aligned}M(u)=&2(b_{0,0}^4-b_{0,0}^1)u+\pi(b_{0,1}^1+a_{1,0}^1)(u^2+u^4)+2(a_{1,0}^1-a_{1,0}^4)u^3\\
&+2(b_{0,1}^4+a_{1,0}^4-b_{0,1}^1-a_{1,0}^1)(u^4+u^2)\int_0^{\frac{1}{\sqrt{1+u^2}}}\sqrt{1-t^2}dt.\end{aligned}\end{eqnarray}
Hence, we get for $h\in(0,\infty)$
$$\frac{d}{du}\Big(\frac{M(u)}{u^2+u^4}\Big)=-\frac{2}{(u+u^3)^2}\big[(b^4_{0,1}-b^1_{0,1})u^4-(a_{1,0}^1-a_{1,0}^4-3b^4_{0,0}+3b^1_{0,0})u^2+
b^4_{0,0}-b^1_{0,0}\big].$$
Thus, by Rolle's theorem, $M(u)$, sa well as $M(h)$, has at most 3 zeros in $(0,+\infty)$. By (3.8), $M(u)$ in (4.6) can be written as
$$\begin{aligned}M(u)=&2(b_{0,0}^4-b_{0,0}^1)u+\frac{\pi}{2}(b_{0,1}^4+a_{1,0}^4+b_{0,1}^1+a_{1,0}^1)(u^2+u^4)+2(a_{1,0}^1-a_{1,0}^4)u^3\\&
-\frac{2}{3}(b_{0,1}^4+a_{1,0}^4-b_{0,1}^1-a_{1,0}^1)u^5+o(u^5)\\
:=&u\big(\tau_0+\tau_1(u+u^3)+\tau_2u^2+\tau_4u^4+o(u^4)\big),\end{aligned}$$
where $\tau_0$, $\tau_1$, $\tau_2$ and $\tau_3$ are constants. It is easy to get that
 $$\det\frac{\partial(\tau_0,\tau_1,\tau_2,\tau_4)}{\partial(b_{0,1}^4,b_{0,1}^1,a_{1,0}^1,b_{0,0}^1)}=\frac{8}{3}\pi.$$
Thus, $\tau_0$, $\tau_1$, $\tau_2$ and $\tau_4$ can be chosen arbitrarily. Therefore, there exist $\tau_0$, $\tau_1$, $\tau_2$ and $\tau_4$ such that $M(u)$ has 3 zeros in $(0,+\infty)$, so does $M(h)$. This completes the proof of Theorem 1.2.\quad $\lozenge$

\vskip 0.2 true cm

\noindent
{\bf Acknowledgment}
 \vskip 0.2 true cm

\noindent
Supported by National Natural Science Foundation of China(11701306,11601250), Construction of First-class Disciplines of Higher Education of Ningxia(Pedagogy)(NXYLXK2017B11) and Ningxia Natural Science Foundation(2019AAC03247) and Key Program of Ningxia Normal University(NXSFZDA1901).

\end{document}